\documentclass[a4paper, fleqn]{article}
\usepackage[english]{babel}
\usepackage{graphicx}
\usepackage{amsmath, amssymb, graphics, float, fullpage}
\usepackage{amsthm}
\usepackage{makeidx}
\usepackage{lipsum}

\newtheorem{theorem}{Theorem}[section]

\newtheorem{proposition}[theorem]{Proposition}
\newtheorem{corollary}[theorem]{Corollary}
\newtheorem{definition}[theorem]{Definition}
\newtheorem{remark}[theorem]{Remark}

\providecommand\phantomsection{}

\makeindex

\begin{document}

\pagenumbering{gobble}
\clearpage
\thispagestyle{empty}

\title{MONODROMY OF THE GENERALIZED HYPERGEOMETRIC EQUATION IN THE\\ FROBENIUS BASIS}
\author{L. D. Molag}
\maketitle
\begin{center}
Department of Mathematics\\
Utrecht University\\
\vspace{1.0cm}

\noindent l.d.molag@uu.nl
\end{center}

\vspace{2cm}
\noindent\textbf{Abstract}\\

\noindent We consider monodromy groups of the generalized hypergeometric equation
\begin{equation*}
\big[z(\theta+\alpha_{1})\cdots (\theta+\alpha_{n})-(\theta+\beta_{1}-1)\cdots (\theta+\beta_{n}-1)\big]f(z) = 0\text{, where }\theta = z d/dz,
\end{equation*}
in a suitable basis, closely related to the Frobenius basis. We pay particular attention to the maximally unipotent case, where $\beta_{1}=\ldots=\beta_{n}=1$, and present a theorem that enables us to determine the form of the corresponding monodromy matrices in the case where $(X-e^{-2\pi i\alpha_{1}})\cdots (X-e^{-2\pi i\alpha_{n}})$ is a product of cyclotomic polynomials.\\

\section{Introduction}
\addcontentsline{toc}{section}{Introduction}

\enlargethispage{\baselineskip}

Let $\alpha_{1},\ldots,\alpha_{n},\beta_{1},\ldots,\beta_{n}\in \mathbb{C}$. The generalized hypergeometric equation
\begin{equation}
\label{generalized hypergeometric equation}
\big[z(\theta+\alpha_{1})\cdots (\theta+\alpha_{n})-(\theta+\beta_{1}-1)\cdots (\theta+\beta_{n}-1)\big]f(z) = 0\text{, where, }\theta = z d/dz
\end{equation}
is a generalization of the Euler-Gauss hypergeometric equation, corresponding to the case $n=2$ which was introduced by Euler in the $18^{th}$ century and studied in the $19^{th}$ century by among others: Gauss, Klein, Riemann and Schwarz.\\

\noindent There exists an $n$-dimensional basis of solutions to (\ref{generalized hypergeometric equation}) in a neighborhood of $z=0$, called the Frobenius basis (at $z=0$). In the case that the local exponents are pairwise distinct (the non-resonant case) this basis is given by $z^{1-\beta_{1}} F_{1},\ldots, z^{1-\beta_{n}} F_{n}$, for some analytic functions $F_{1},\ldots,F_{n}$, known as Clausen-Thomae hypergeometric functions, that are defined on some open neighborhood of $0$. In the case that all local exponents equal $1$ (the maximally unipotent case) the Frobenius basis is of the following form:
\begin{align*}
f_{0} &= 1+h_{0}\\
f_{1} &= f_{0}\log(z)+ h_{1}\\
f_{2} &= \frac{1}{2}  f_{0} \log^{2}(z)+h_{1} \log(z)+h_{2}\\
& \vdots\\
f_{n-1} &=  \frac{1}{(n-1)!} f_{0} \log^{n-1}(z)+ \sum_{l=0}^{n-2} \frac{1}{l!} h_{n-1-l} \log^{l}(z).
\end{align*}
where the $h_{l}$ are analytic, vanishing in $z=0$, and the unique functions with this property.\\

\noindent We are mainly interested in the monodromy corresponding to the Frobenius basis. Important to us will be the explicit form of matrices that are used in the proof of Levelt's theorem\cite{unpublished}, from which one can deduce the explicit form of the monodromy matrices corresponding to (\ref{generalized hypergeometric equation}) in a certain basis. It turns out that we can actually find the corresponding basis of functions explicitly, these functions are known as Melllin-Barnes integrals and the corresponding basis is called the Mellin-Barnes basis. The advantage of this basis is that the functions are defined on a large region, whereas the functions in the Frobenius basis are generally determined by powerseries with finite convergence radius (although they can be analytically extended). Our intention of course, is to express the functions in the Frobenius basis as linear combinations of Mellin-Barnes integrals, such that we can easily continue them along a path. In the next chapter it will be explained in detail how this is done. \\  

\noindent In the non-resonant case it follows immediately that the monodromy matrix  around $0$ in the Frobenius basis around $z=0$ equals $\text{diag}(e^{-2\pi i\beta_{1}},\ldots,e^{-2\pi i\beta_{n}})$. Theorem \ref{nonresonantthm} explains the general structure of the monodromy group, by giving the explicit form of the monodromy matrix around $1$ in the Frobenius basis around $z=0$, namely its $(k,l)$ entry, with $k,l=1,2,\ldots,n$, is
\begin{align}
\delta_{kl}+c e^{2\pi i\beta_{k}} \prod_{m=1}^{n} \frac{\sin(\pi(\beta_{l}-\alpha_{m}))}{\sin(\pi(\beta_{l}-\beta_{m}))}.
\end{align}
Here $c=2i (-1)^{n} e^{\pi i(\beta_{1}-\alpha_{1}+\ldots+\beta_{n}-\alpha_{n})}$ and the factor $\sin(\pi(\beta_{l}-\beta_{l}))$ should be read as $1$. This shows in particular that all monodromy matrices have algebraic entries when the parameters $\alpha_{1},\ldots,\alpha_{n},\beta_{1},\ldots,\beta_{n}$ are rational, a property that is not shared with the maximally unipotent case.\\

\noindent Our main theorem, about the maximally unipotent case, Theorem \ref{main?}, will need the following result. Suppose that $\alpha_{1},\ldots,\alpha_{n}\in\mathbb{C}\setminus\mathbb{Z}$ are such that $(X-e^{-2\pi i\alpha_{1}})\cdots (X-e^{-2\pi i\alpha_{n}})$ is a product of cyclotomic polynomials, then we can find a number $r\in\mathbb{N}$ and numbers $a_{1},\ldots,a_{r},b_{1},\ldots,b_{r}\in\mathbb{N}$ such that
\begin{align*}
(X-e^{-2\pi i\alpha_{1}})\cdots (X-e^{-2\pi i\alpha_{n}}) = \frac{X^{a_{1}}-1}{X^{b_{1}}-1}\cdots \frac{X^{a_{r}}-1}{X^{b_{r}}-1}.
\end{align*}
When this is the case it will turn out that, equivalently, we could investigate the equation
\begin{align*}
\theta^{n}f=Cz(\theta-\alpha_{1})\cdots (\theta-\alpha_{n})f\text{ where }C = \frac{a_{1}^{a_{1}}\cdots a_{n}^{a_{n}}}{b_{1}^{b_{1}}\cdots b_{n}^{b_{n}}},
\end{align*}
which has its own Frobenius basis $f_{n-1}^{C},\ldots,f_{1}^{C},f_{0}^{C}$. This corresponds to the normalization $z\to Cz$, i.e. $f_{k}^{C}(z)=f_{k}(C z)$ for $k=0,\ldots,n-1$.   In fact this is precisely what the authors of \cite{CYY} do for the case $n=4$, in that case the hypergeometric equations arise from Calabi-Yau threefolds. They showed, using a basis that shows resemblance to the Mellin-Barnes basis, that the entries of the corresponding monodromy matrices contain geometric invariants of these Calabi-Yau threefolds. In particular, they gave a neat expression for the monodromy matrices. Generalization of their result for arbitrary $n$ has been our motivation to study the maximally unipotent case.

Our main theorem gives us insight in to the general form of the monodromy matrices in the case that $(z-e^{-2\pi i\alpha_{1}})\cdots (z-e^{-2\pi i\alpha_{n}})$ defines a product of cyclotomic polynomials, in particular it provides us with a practical method to determine the monodromy matrices. We will see that all matrices in the corresponding monodromy group have their entries in $\mathbb{Q}(\zeta(3) (2\pi i)^{-3},\zeta(5) (2\pi i)^{-5},\ldots,\zeta(m)(2\pi i)^{-m})$, with $m$ the largest odd number below $n$.\\

I would like to thank Frits Beukers, who was the supervisor of my master thesis, which contains a lot of material that is being used in this article, for advising me to publish my results and helping me along the way. I am thankful to Willem Pranger for pointing out numerous issues for substantive improvement in my master thesis, and consequently this article. I thank Julian Lyczak and Merlijn Staps for their proof of theorem \ref{integerC}.

\section{Monodromy groups of the generalized hypergeometric equation}

\subsection{The Mellin-Barnes basis}

Let $z_{0}$ be an element of $\{0,1,\infty\}$, the set of singularities corresponding to (\ref{generalized hypergeometric equation}). We will denote the monodromy matrix around $z_{0}$ by $M_{z_{0}}$. For (\ref{generalized hypergeometric equation}) we know that $M_{0}$ has eigenvalues $e^{-2\pi i\beta_{1}},\ldots,e^{-2\pi i\beta_{n}}$ and $M_{\infty}$ has eigenvalues $e^{2\pi i\alpha_{1}},\ldots,e^{2\pi i\alpha_{n}}$. We will consider the case where all eigenvalues  $e^{-2\pi i\beta_{1}},\ldots,e^{-2\pi i\beta_{n}}$ differ from the eigenvalues $e^{2\pi i\alpha_{1}},\ldots,e^{2\pi i\alpha_{n}}$. Here and in the rest of this article we will demand that these two sets of eigenvalues are disjunct, i.e. $\alpha_{k}$ differs from $\beta_{l}$ modulo $1$ for all $k,l=1,2,\ldots,n$. A matrix will be called a (pseudo-)reflection if this matrix minus the identity has rank $1$. The following theorem gives us insight in to the general form of the monodromy matrices corresponding to this case. 

\begin{theorem}
\label{Levelt}
\textbf{(Levelt)}\index{Levelt's theorem} Let $a_{1},\ldots,a_{n},b_{1},\ldots,b_{n}\in\mathbb{C}\setminus\{0\}$ be such that $a_{i}\neq b_{j}$ for all $1\leq i,j\leq n$. Then there exist $A,B\in GL(n,\mathbb{C})$ with eigenvalues $a_{1},\ldots,a_{n}$ and $b_{1},\ldots,b_{n}$ respectively such that $AB^{-1}$ is a reflection. Moreover, the pair $A, B$ is uniquely determined up to conjugation. 
\end{theorem}

\noindent What is important about Levelt's theorem is its proof \cite{unpublished}. It shows us explicitly what the monodromy matrices look like in a particular basis chosen, namely
\begin{align*}
A = \left(\begin{array}{ccccc} 0 & 1 & 0 & \ldots & 0\\
0 & 0 & 1 & \ldots & 0\\
\vdots & & &  & \vdots\\
0 & 0 & 0 & \ldots & 1\\
-A_{n} & -A_{n-1} & -A_{n-2} & \ldots & -A_{1}
\end{array}\right)
\text{ and }
B = \left(\begin{array}{ccccc} 0 & 1 & 0 & \ldots & 0\\
0 & 0 & 1 & \ldots & 0\\
\vdots & & &  & \vdots\\
0 & 0 & 0 & \ldots & 1\\
-B_{n} & -B_{n-1} & -B_{n-2} & \ldots & -B_{1}
\end{array}\right),
\end{align*}
where $A_{1},\ldots, A_{n}, B_{1},\ldots, B_{n}$ are defined through $(X-a_{1})\cdots (X-a_{n})=X^{n}+A_{1} X^{n-1}+\ldots+A_{n}$ and $(X-b_{1})\cdots (X-b_{n})=X^{n}+B_{1} X^{n-1}+\ldots+B_{n}$.\\

\noindent It is known that $M_{1}$ has $n-1$ eigenvalues equal to $1$ and is thus a reflection (and so is $M_{\infty}^{-1} M_{1}^{-1} M_{\infty}$). In particular $M_{0}$ and $M_{\infty}^{-1}$, satisfying the relation $M_{0}M_{1}M_{\infty}=\mathbb{I}$, play the role of $A$ and $B$ in Levelt's theorem. It turns out that we can actually find an explicit basis of functions in which $M_{0}$ equals the matrix $A$ used in Levelt's theorem, with $a_{k}=e^{-2\pi i\beta_{k}}$ for $k=1,\ldots,n$. In the following we will choose the argument of $z$ in $(0,2\pi)$, which determines $z^{s}=|z|^{s} e^{i\text{arg}(z) s}$. 

\begin{definition}
Let $\alpha_{1},\ldots,\alpha_{n},\beta_{1},\ldots,\beta_{n}\in\mathbb{C}$ and $\alpha_{k}$ differs from $\beta_{l}$ modulo $1$ for all $k,l=1,2,\ldots,n$. We define for $j=0,1,\ldots,n-1$ and $z\in\mathbb{C}\setminus\mathbb{R}_{\geq 0}$
\begin{equation}
\label{MellinBarnes}
I_{j}(z) = \frac{(-1)^{n}}{(2\pi i)^{n}}\int_{L} \left(\prod_{k=1}^{n} \Gamma(\alpha_{k}+s)\Gamma(1-\beta_{k}-s)\right) e^{(2j-n)\pi i s}z^{s} ds.
\end{equation}
Here $L$ is a path from $i\infty$ to $-i\infty$ that bends in such a way that all points $-\alpha_{k}-m$ with $m\in\mathbb{Z}_{\geq 0}$ are on the left of it and all points $1-\beta_{k}+m$ with $m\in\mathbb{Z}_{\geq 0}$ are on the right of it, for big enough $s$ we require it to be on the imaginary axis.\\
\end{definition}

\begin{remark}
Here by `left' and `right' we mean that $L$ divides $\mathbb{C}\setminus L$ into two connected components, the component that contains all $s$ with negative real part for $s$ big enough will be referred to as the left component, the other as the right component. The requirement that $L$ is on the imaginary axis for big $s$ is not necessary but will turn out to be convenient in what follows. 
\end{remark}

\noindent Let us argue that the Mellin-Barnes integrals (\ref{MellinBarnes}) are well defined. Stirling's formula\index{Stirling's formula} tells us that for $a,b\in\mathbb{R}$, $a$ bounded, we have
\begin{align*}
|\Gamma(a+bi)| = \mathcal{O}(|b|^{a-1/2} e^{-\pi |b|/2})\text{ as }|b|\to\infty.
\end{align*}
We deduce that $|\Gamma(\alpha_{k}+it)\Gamma(1-\beta_{k}-it)|=\mathcal{O}(|t|^{1+\Re(\alpha_{k}-\beta_{k})} e^{-\pi|t|})$ as $|t|\to\infty$. Henceforth for $j=0,1,\ldots,n-1$
\begin{align}
\label{stirling}
\left| \left(\prod_{k=1}^{n} \Gamma(\alpha_{k}+it)\Gamma(1-\beta_{k}-it)\right) e^{((2j-n)\pi i)it}(i t)^{s}\right|\\
= \mathcal{O}(|t|^{n+\sum_{k=1}^{n} \Re(\alpha_{k}-\beta_{k})} e^{-\text{arg}(z)\pi |t|})\text{ as }|t|\to\infty.
\end{align}
Since the argument of $z$ is positive we conclude that the integrals $I_{j}$ converge.\\

\begin{proposition}
\label{MBversch}
Let $N\in\mathbb{N}$. Denote by $i_{j,z}$ the integrant of $I_{j}(z)$. Define by $R(N)$ the set of singularities of $i_{j,z}(s)$ between $L$ and $L+N$ and by $R(\infty)$ and $R(-\infty)$ the set of singularities on the right respectively on the left of $L$. Denote by $I_{j}^{N}$ the integral $I_{j}$ were the path $L$ has been replaced by $L+N$. We have (for a fixed choice of $\pm$)
\begin{align}
I_{j}(z) = I_{j}^{\pm N}(z) \pm 2\pi i\sum_{p\in R(\pm N)} \text{Res}_{p}( i_{j,z}),
\end{align}
In particular we have for $|z|^{\pm 1}<1$ that
\begin{align}
I_{j}(z) =\pm 2\pi i\sum_{p\in R(\pm \infty)} \text{Res}_{p}(i_{j,z}).
\end{align}
\end{proposition}

\noindent\textbf{Proof. }For $T>0$ big enough consider the path $L(T)$ that coincides with $L$ but is from $iT$ to $-iT$. Now connect the paths $L(T)$ and $L(T)\pm N$ (for a fixed choice of $\pm$) by two linear segments $L_{-}(T)$ and $L_{+}(T)$ from $-iT$ to $\pm N-iT$ and from $\pm N+iT$ to $iT$ respectively. Thus we get a closed path and by the residue theorem
\begin{align*}
\int_{L(T)}+\int_{L_{-}(T)}-\int_{L(T)\pm N}+\int_{L_{+}(T)} i_{j,z}(s) ds = \pm 2\pi i\sum_{p\in R(\pm N)} \text{Res}_{p}( i_{j,z}). 
\end{align*}
For the first part of the proposition it suffices to show that the integrals over $L_{\pm}(T)$ tend to $0$ as $T\to\infty$. For this we use the Stirling approximation: $|i_{j,z}(t\pm iT)|= \mathcal{O}(T^{n+2nN+\sum_{k=1}^{n} \Re(\alpha_{k}-\beta_{k})} e^{-\text{arg}(z)\pi T})$. This tends to $0$ as $T\to\infty$, as the integration intervals are finite this proves that the integrals over $L_{\pm}(T)$ tend to $0$ as $T\to\infty$.\\

\noindent Now for the second part of the proposition we should prove that the integral over $L\pm N$ tends to $0$ as $N\to\infty$ whenever $|z|^{\pm 1}<1$. We will prove this only for the $|z|<1$ case, the other case is analogous. We see that for $s$ on $L$ we have
\begin{align*}
|\Gamma(\alpha_{k}+s+N)\Gamma(1-\beta_{k}-(s+N))| = \left|\prod_{j=0}^{N-1}\frac{\alpha_{k}+s+j}{1-\beta_{k}+s+j}\right| |\Gamma(\alpha_{k}+s)\Gamma(1-\beta_{k}-s)|.
\end{align*}
We notice that uniformly on $L$
\begin{align*}
\lim_{j\to\infty} \left|\frac{\alpha_{k}+s+j}{1-\beta_{k}+s+j}\right|\leq \lim_{j\to\infty} 1+\frac{|\alpha_{k}+1-\beta_{k}|}{|1-\beta_{k}+s+j|} = 1,
\end{align*}
where we have used that the real part of $s$ is bounded on $L$. In particular for $j$ big enough we have uniformly on $L$ that
\begin{align*}
\left|\frac{\alpha_{k}+s+j}{1-\beta_{k}+s+j}\right|\leq |z|^{-\frac{1}{2n}}.
\end{align*}
We conclude that the integrant of the integral over $L+N$ satisfies the same inequality as in (\ref{stirling}), but with a factor $|z|^{\frac{N}{2}}$ in front of it. Since $|z|<1$ we conclude that the integral over $L+N$ converges to $0$. 
\begin{flushright}$\square$\end{flushright}

\begin{theorem}
\label{MBbasis}
The functions $I_{0},\ldots,I_{n-1}$ form a basis $\mathcal{I}$, the Mellin-Barnes basis\index{Mellin-Barnes basis}, of the generalized hypergeometric equation (\ref{generalized hypergeometric equation}).
\end{theorem}

\noindent\textbf{Proof. }Let us prove that they are solutions to the generalized hypergeometric equation. First we notice that
\begin{align*}
\theta e^{(2j-n)\pi i s}z^{s} = z  e^{(2j-n)\pi i s} s z^{s-1} = s e^{(2j-n)\pi i s} z^{s}.
\end{align*}
Thus
\begin{align*}
z(\theta+\alpha_{1})\cdots (\theta+\alpha_{n}) I_{j} &=  \frac{(-1)^{n}}{(2\pi i)^{n}}\int_{L} \left(\prod_{k=1}^{n} \Gamma(\alpha_{k}+s)\Gamma(1-\beta_{k}-s)\right)\\
&\times (s+\alpha_{1})\cdots (s+\alpha_{n}) e^{(2j-n)\pi i s} z^{s+1} ds
\end{align*}
\begin{align*}
&=  \frac{(-1)^{n}}{(2\pi i)^{n}}\int_{L} \left(\prod_{k=1}^{n}\Gamma(\alpha_{k}+s+1)\Gamma(1-\beta_{k}-s)\right)  e^{(2j-n)\pi i s} z^{s+1} ds\\
&= \frac{(-1)^{n}}{(2\pi i)^{n}}\int_{L+1} \left(\prod_{k=1}^{n} \Gamma(\alpha_{k}+s)(1-\beta_{k}-s)\Gamma(1-\beta_{k}-s)\right) e^{(2j-n)\pi i s}(-1)^{n} z^{s} ds\\
&= (\theta+\beta_{1}-1)\cdots (\theta+\beta_{n}-1)I_{j}^{1}(z)\\
&+ 2\pi i (\theta+\beta_{1}-1)\cdots (\theta+\beta_{n}-1)\sum_{p\in R(1)} \text{Res}_{p}( i_{j,z})
\end{align*}
by Proposition \ref{MBversch}. Now if there are indeed singularities in $R(1)$ they must be of the form $s=1-\beta_{k}$. The Residue corresponding to such a pole is a linear combination of terms of the form $\log^{l}(z) z^{1-\beta_{k}}$ for $0\leq l<n$. If such a term appears then $\beta_{k}$ must have degeneracy at least $l+1$. We notice using the Leibniz rule that
\begin{align*}
& (\theta+\beta_{k}-1)^{l+1} \log^{l}(z) z^{1-\beta_{k}} =  (\theta+\beta_{k}-1)^{l} l\log^{l-1}(z) z^{1-\beta_{k}}\\
&= \ldots = (\theta+\beta_{k}-1) l(l-1)\cdots 1\cdot z^{1-\beta_{k}}=0.
\end{align*}
Hence
\begin{align*}
&(\theta+\beta_{1}-1)\cdots (\theta+\beta_{n}-1)\\
&\frac{(-1)^{n}}{(2\pi i)^{n}} \int_{L+1} \left(\prod_{k=1}^{n} \Gamma(\alpha_{k}+s)\Gamma(1-\beta_{k}-s)\right) e^{(\log(z)+(2j-n)\pi i)s} ds\\
&=  (\theta+\beta_{1}-1)\cdots (\theta+\beta_{n}-1) I_{j}(z)
\end{align*}
and we conclude that the $I_{j}$ are solutions to the hypergeometric equation. Suppose $I_{0},\dots,I_{n-1}$ do not form a basis. Then there exists a polynomial $p$ of degree at most $n-1$, not identically zero, such that
\begin{align*}
\int_{L} \left(\prod_{k=1}^{n} \Gamma(\alpha_{k}+s)\Gamma(1-\beta_{k}-s)\right) e^{-\pi i n s}p(e^{2\pi i s}) z^{s} ds = 0.
\end{align*}
This is only possible if no terms of the form $\log^{l}(z) z^{1-\beta_{k}}$ occur (when evaluated in a neighborhood of $z=0$), i.e. that all singularities of the original integrant are removed by $p(e^{2\pi i s})$ (see remark \ref{combilogremark} for clarification). This implies that $p$ must have all $e^{-2\pi i\beta_{k}}$ as roots (with the same multiplicity as $\beta_{k}$), and this is a contradiction since it requires $p$ to have degree at least $n$. 
\begin{flushright}$\square$\end{flushright}

\begin{theorem}
\label{MBmonodromy}
Suppose $\alpha_{k}$ differs from the $\beta_{l}$ modulo $1$ for all $1\leq k,l\leq n$. The monodromy matrices in the Mellin-Barnes basis are
\begin{align*}
M_{0} &= \left( \begin{array}{ccccc}
0 & 1 & 0 & \hdots & 0\\
0 & 0 & 1 & \hdots & 0\\
\vdots & \vdots &  &  & \vdots\\
0 & 0 & 0 & \hdots & 1\\
-B_{n} & -B_{n-1} & -B_{n-2} & \hdots & -B_{1} \end{array} \right)
\end{align*}
\begin{align*}
M_{1} &= \left( \begin{array}{ccccc}
1+\frac{A_{n}-B_{n}}{B_{n}} & \frac{A_{n-1}-B_{n-1}}{B_{n}} & \frac{A_{n-2}-B_{n-2}}{B_{n}} & \hdots & \frac{A_{1}-B_{1}}{B_{n}}\\
0 & 1 & 0 & \hdots & 0\\
0 & 0 & 1 & \hdots & 0\\
\vdots & \vdots &  &  & \vdots\\
0 & 0 & 0 & \hdots & 1 \end{array} \right)
\end{align*}
\begin{align*}
M_{\infty} &= \left( \begin{array}{ccccc}
-\frac{A_{n-1}}{A_{n}} & -\frac{A_{n-2}}{A_{n}} & -\frac{A_{n-3}}{A_{n}} & \hdots & -\frac{A_{0}}{A_{n}}\\
1 & 0 & 0 & \hdots & 0\\
0 & 1 & 0 & \hdots & 0\\
\vdots & \vdots &  &  & \vdots\\
0 & 0 & \hdots & 1 & 0 \end{array} \right)
\end{align*}
Where $z^{n}+B_{1}z^{n-1}+\ldots+B_{n-1}z+B_{n}$ is the polynomial with roots $e^{-2\pi i\beta_{k}}$, $k=1,2,\ldots,n$ and $z^{n}+A_{1}z^{n-1}+\ldots+A_{n-1}z+A_{n}$ is the polynomial with roots $e^{-2\pi i\alpha_{k}}$, $k=1,2,\ldots,n$.\\
\end{theorem}

\noindent\textbf{Proof. }By construction we have $I_{j}\to I_{j+1}$ under a counterclockwise loop around $0$ for $j=0,1,\ldots,n-2$. Notice that
\begin{align*}
-B_{n}I_{0}-\ldots-B_{1} I_{n-1} &= \frac{(-1)^{n}}{(2\pi i)^{n}}\int_{L} \left(\prod_{k=1}^{n} \Gamma(\alpha_{k}+s)\Gamma(1-\beta_{k}-s)\right) e^{-\pi i n s} z^{s}\\
& \times \left(e^{2\pi i n s} - \prod_{k=1}^{n} (e^{2\pi i s}-e^{-2\pi i\beta_{k}})\right) ds.
\end{align*}
Notice what happens when we lower the argument by $2\pi$. By the same arguments used in the proof of Proposition \ref{MBversch} we have that
\begin{align*}
\frac{(-1)^{n}}{(2\pi i)^{n}}\int_{L} \left(\prod_{k=1}^{n}\Gamma(\alpha_{k}+s)\Gamma(1-\beta_{k}-s)\right) e^{-2\pi i s-\pi i n s} z^{s} \prod_{k=1}^{n} (e^{2\pi i s}-e^{-2\pi i\beta_{k}}) ds
\end{align*}
is equal to $2\pi i$ times the sum of its residues corresponding to its singularities to the right of $L$ for $|z|<1$. But it has no (non removable) singularities in that region so it vanishes. We conclude that when we lower the argument by $2\pi$ then $-B_{n}I_{0}-\ldots-B_{1} I_{n-1}$ transforms to $I_{n-1}$, i.e. a counterclockwise loop around the origin corresponds to the transformation $I_{n-1}\to -B_{n}I_{0}-\ldots-B_{1} I_{n-1}$.\\ 

\noindent From the Frobenius basis around $\infty$ it is clear that $M_{\infty}^{-1}$ should have eigenvalues $e^{-2\pi i \alpha_{1}},\ldots,e^{-2\pi i \alpha_{n}}$. Furthermore, we know that  $M_{0} M_{\infty} = M_{\infty}^{-1} M_{1}^{-1}  M_{\infty}$ is a reflection. Hence we may apply Levelts theorem (\ref{Levelt}) to conclude that
\begin{align*}
(M_{\infty})^{-1} =\left( \begin{array}{ccccc}
0 & 1 & 0 & \hdots & 0\\
0 & 0 & 1 & \hdots & 0\\
\vdots & \vdots &  &  & \vdots\\
0 & 0 & 0 & \hdots & 1\\
-A_{n} & -A_{n-1} & -A_{n-2} & \hdots & -A_{1} \end{array} \right)
\end{align*}
The forms of $M_{\infty}$ and $M_{1}$ now easily follow. 
\begin{flushright}$\square$\end{flushright}

\subsection{The non-resonant case}

In this section we will consider the case where $\beta_{1},\ldots,\beta_{n}$ are distinct modulo $1$ and the $\alpha_{1},\ldots,\alpha_{n}$ are distinct from the $\beta_{1},\ldots,\beta_{n}$ modulo $1$\index{Non-resonant}. Though our research is mainly aimed at the maximally unipotent case, we treat the non-resonant case because it is barely any extra work, and the results can be compared with that of the maximally unipotent case. In the Frobenius basis at $0$, denoted by $f_{1},\ldots,f_{n}$, we have
\begin{align*}
M_{0} = \left(\begin{array}{cccc} 
e^{-2\pi i\beta_{1}} & 0 & \hdots & 0\\
0 & e^{-2\pi i\beta_{2}} & \hdots & 0\\
\vdots & & \ddots & 0\\
0 & 0 & \hdots & e^{-2\pi i\beta_{n}}\end{array}\right).
\end{align*}
We would also like to express the monodromy matrices $M_{1}$ and $M_{\infty}$ in the Frobenius basis at $z=0$. For this purpose we will prove the following theorem about the transformation matrix\index{Transformation matrix} between the Mellin-Barnes basis and the Frobenius basis at $z=0$. 

\begin{proposition}
We have
\begin{align}
\left(\begin{array}{c}I_{0}\\ \vdots \\ I_{n-1}\end{array}\right) = V D \left(\begin{array}{c}f_{1}\\ \vdots \\ f_{n}\end{array}\right)
\end{align}
where $V$ is the VanderMonde matrix\index{VanderMonde matrix} $V_{kl} = e^{-2\pi i k\beta_{l}}$ and $D$ is the diagonal matrix with entries
\begin{align*}
D_{ll} = \frac{1}{(2i)^{n-1}} e^{\pi i(n-2k)\beta_{l}} \frac{\Gamma(\alpha_{1}-\beta_{l}+1)\cdots \Gamma(\alpha_{n}-\beta_{l}+1)}{\Gamma(\beta_{1}-\beta_{l}+1)\cdots \Gamma(\beta_{n}-\beta_{l}+1)} \left(\prod_{m=1,m\neq l}^{n} \frac{1}{\sin(\pi(\beta_{m}-\beta_{l}))}\right)
\end{align*}
with $k=0,1,\ldots,n-1$ and $l=1,\ldots,n$. 
\end{proposition}

\noindent\textbf{Proof. }Using Proposition \ref{MBversch} we conclude that
\begin{align*}
I_{k} &= \frac{(-1)^{n}}{(2\pi i)^{n-1}} \sum_{l=1}^{n} \sum_{m=0}^{\infty} \lim_{s\to 1-\beta_{l}+m} (s-1+\beta_{l}-m)\Gamma(1-\beta_{l}-s)\\
&\times \Gamma(\alpha_{l}+s)\left(\prod_{p=1,p\neq l}^{n} \Gamma(\alpha_{p}+s)\Gamma(1-\beta_{p}-s)\right) e^{(2k-n)\pi i s}z^{s}\\
&= \frac{1}{(2\pi i)^{n-1}} \sum_{l=1}^{n} \sum_{m=0}^{\infty} \frac{(-1)^{m}}{m!} e^{\pi i(n-2k)\beta_{l}} (-1)^{nm} z^{1-\beta_{l}+m}\\
&\times  \Gamma(\alpha_{l}-\beta_{l}+1+m)\left(\prod_{p=1,p\neq l}^{n}\Gamma(\alpha_{p}-\beta_{l}+1+m)\Gamma(\beta_{l}-\beta_{p}-m)\right)\\
&= \frac{1}{(2i)^{n-1}} \sum_{l=1}^{n} e^{\pi i(n-2k)\beta_{l}} \left(\prod_{p=1,p\neq l}^{n} \frac{1}{\sin(\pi(\beta_{p}-\beta_{l}))}\right)\\
&\times \frac{\Gamma(\alpha_{1}-\beta_{l}+1)\cdots \Gamma(\alpha_{n}-\beta_{l}+1)}{\Gamma(\beta_{1}-\beta_{l}+1)\cdots \Gamma(\beta_{n}-\beta_{l}+1)} z^{1-\beta_{l}}\sum_{m=0}^{\infty} \frac{(\alpha_{1}-\beta_{l}+1)_{m}\cdots (\alpha_{n}-\beta_{l}+1)_{m}}{(\beta_{1}-\beta_{l}+1)_{m}\cdots (\beta_{n}-\beta_{l}+1)_{m}} z^{m}.
\end{align*}
Therefore
\begin{align*}
I_{k} &= \frac{1}{(2i)^{n-1}} \sum_{l=1}^{n} e^{\pi i(n-2k)\beta_{l}} \frac{\Gamma(\alpha_{1}-\beta_{l}+1)\cdots \Gamma(\alpha_{n}-\beta_{l}+1)}{\Gamma(\beta_{1}-\beta_{l}+1)\cdots \Gamma(\beta_{n}-\beta_{l}+1)} \left(\prod_{p=1,p\neq l}^{n} \frac{1}{\sin(\pi(\beta_{p}-\beta_{l}))}\right) f_{l}(z).
\end{align*}
\begin{flushright}$\square$\end{flushright}

\begin{theorem}
\label{nonresonantthm}
Define $c=2i (-1)^{n} e^{\pi i(\beta_{1}-\alpha_{1}+\ldots+\beta_{n}-\alpha_{n})}$. In the Frobenius basis at $z=0$ the monodromy matrix around $z=1$ satisfies
\begin{align}
(M_{1})_{kl} &= \delta_{kl}+c e^{2\pi i\beta_{k}} \prod_{m=1}^{n} \frac{\sin(\pi(\beta_{l}-\alpha_{m}))}{\sin(\pi(\beta_{l}-\beta_{m}))}
\end{align}
where $k,l=1,2,\ldots,n$ and the factor $\sin(\pi(\beta_{l}-\beta_{l}))$ should be read as $1$. 
\end{theorem}

\noindent\textbf{Proof. }We calculate
\begin{align*}
\sum_{m=0}^{n-1} \frac{A_{n-m}-B_{n-m}}{B_{n}} e^{-2\pi i m\beta_{l}} &= \frac{1}{B_{n}} \left(\prod_{m=1}^{n} (e^{-2\pi i\beta_{l}}-e^{-2\pi i\alpha_{m}})-\prod_{m=1}^{n} (e^{-2\pi i\beta_{l}}-e^{-2\pi i\beta_{m}})\right)\\
&= (2i)^{n}e^{2\pi i(\beta_{1}+\ldots+\beta_{n})} e^{-\pi i(\alpha_{1}+\ldots+\alpha_{n})} e^{-\pi i n \beta_{l}} \prod_{m=1}^{n} \sin(\pi(\alpha_{m}-\beta_{l}))\\
&= 2i e^{2\pi i(\beta_{1}+\ldots+\beta_{n})} e^{-\pi i(\alpha_{1}+\ldots+\alpha_{n})} \tilde{D}^{-1}_{ll} \sin(\pi(\alpha_{l}-\beta_{l})) \prod_{m=1,m\neq l}^{n} \frac{\sin(\pi(\alpha_{m}-\beta_{l}))}{\sin(\pi(\beta_{m}-\beta_{l}))}.
\end{align*}
where
\begin{align*}
\tilde{D}_{ll}\frac{\Gamma(\alpha_{1}-\beta_{l}+1)\cdots \Gamma(\alpha_{n}-\beta_{l}+1)}{\Gamma(\beta_{1}-\beta_{l}+1)\cdots \Gamma(\beta_{n}-\beta_{l}+1)} = D_{ll}.
\end{align*}
To complete the proof we will have to determine the inverse of $V$. We notice that this inverse is determined by
\begin{align*}
\prod_{m=1,m\neq k}^{n} \frac{z-e^{-2\pi i\beta_{m}}}{e^{-2\pi i\beta_{k}}-e^{-2\pi i\beta_{m}}} = (V^{-1})_{k,0}+(V^{-1})_{k,1} z+\ldots+(V^{-1})_{k,n-1} z^{n-1}.
\end{align*}
We will only need the first column of $V^{-1}$, the $k^{th}$ entry of this column is
\begin{align*}
\prod_{m=1,m\neq k}^{n} \frac{-e^{-2\pi i\beta_{m}}}{e^{-2\pi i\beta_{k}}-e^{-2\pi i\beta_{m}}} = (-1)^{n-1}e^{2\pi i\beta_{k}} e^{-\pi i(\beta_{1}+\ldots+\beta_{n})} \tilde{D}_{kk}.
\end{align*}
We conclude that the matrix $M_{1}-\mathbb{I}$ equals $2i (-1)^{n-1}e^{\pi(\beta_{1}-\alpha_{1}+\ldots+\beta_{n}-\alpha_{n})}$ times
\begin{align}
\left(\begin{array}{c} e^{2\pi i\beta_{1}} \tilde{D}_{11}\\ \vdots \\ e^{2\pi i\beta_{n}} \tilde{D}_{nn} \end{array}\right)
\left(\begin{array}{c} \sin(\pi(\alpha_{1}-\beta_{1})) \tilde{D}_{11}^{-1} \prod_{m=1,m\neq 1}^{n} \frac{\sin(\pi(\alpha_{m}-\beta_{1}))}{\sin(\pi(\beta_{m}-\beta_{1}))} \\ \vdots \\ \sin(\pi(\alpha_{n}-\beta_{n})) \tilde{D}_{nn}^{-1} \prod_{m=1,m\neq n}^{n} \frac{\sin(\pi(\alpha_{m}-\beta_{n}))}{\sin(\pi(\beta_{m}-\beta_{n}))}\end{array}\right)^{T}
\end{align}
which implies the desired result. 
\begin{flushright}$\square$\end{flushright}

Though the form of $M_{1}$ is the easiest to find the following proposition will show that the form of $M_{\infty}$ can easily be deduced from the form of $M_{1}$. 

\begin{proposition}
\label{inverserank1}
Let $M$ be an $n\times n$ matrix with rank $\leq 1$. Suppose that $\mathbb{I}+M$ is invertible. Then
\begin{align*}
(\mathbb{I}+M)^{-1} = \mathbb{I}-\frac{1}{1+\text{Tr}(M)}M.
\end{align*}
\end{proposition}

\noindent\textbf{Proof. }Since $M$ has rank $\leq 1$ it can be written as $M_{kl}=u_{k}v_{l}$ for $n$-dimensional vectors $u$ and $v$. Thus we notice that
\begin{align*}
(M^{2})_{kl} = \sum_{m=1}^{n} u_{k}v_{m}u_{m}v_{l} = \text{Tr}(M) M_{kl}.
\end{align*}
Since $M$ has rank $\leq 1$ we know that it has $n-1$ eigenvalues equal to $0$. The condition that $\mathbb{I}+M$ is invertible thus boils down to $\text{Tr}(M)\neq -1$. We see that
\begin{align*}
(\mathbb{I}+M)(\mathbb{I}-\frac{1}{1+\text{Tr}(M)}M) = \mathbb{I}+M-\frac{1}{1+\text{Tr}(M)} M-\frac{\text{Tr}(M)}{1+\text{Tr}(M)}M=\mathbb{I}.
\end{align*}
\begin{flushright}$\square$\end{flushright}

\begin{corollary}
Suppose $\alpha_{1},\ldots,\alpha_{n},\beta_{1},\ldots,\beta_{n}$ are distinct modulo $1$. Then in the Frobenius basis at $z=0$ the monodromy matrix around $z=\infty$ satisfies
\begin{align}
(M_{\infty})_{kl} = e^{2\pi i\alpha_{k}}\delta_{kl}+\frac{4}{c} e^{2\pi i(\beta_{k}+\alpha_{k})} \prod_{m=1}^{n} \frac{\sin(\pi(\beta_{l}-\alpha_{m}))}{\sin(\pi(\beta_{l}-\beta_{m}))}
\end{align}
where $k,l=1,2,\ldots,n$. 
\end{corollary}

\noindent\textbf{Proof. }We know that $1+\text{Tr}(M_{1}-\mathbb{I})=1+(A_{n}-B_{n})/B_{n} = -c^{2}/4$. Hence
\begin{align*}
M_{\infty} &= (\mathbb{I}+4 c^{-2}(M_{1}-\mathbb{I})) M_{0}^{-1},
\end{align*}
leading to the desired result. 
\begin{flushright}$\square$\end{flushright}

We conclude this paragraph with the remark that when $\alpha_{1},\ldots,\alpha_{n},\beta_{1},\ldots,\beta_{n}\in\mathbb{Q}$ the corresponding monodromy group consists of matrices with algebraic entries. In the next chapter it will become clear that this is no longer implied in the maximally unipotent case. 

\subsection{The maximally unipotent case}

In this section we will consider the case where $\beta_{1}=\ldots=\beta_{n}=1$\index{Maximally unipotent}. \noindent In what follows it will turn out that our results become more elegant when we slightly alter the Frobenius basis. We will consider the ordered basis $\{f_{n-1}/(2\pi i)^{n-1},f_{n-2}/(2\pi i)^{n-2},\ldots,f_{0}\}$ instead. Notice that in this basis we have
\begin{align*}
M_{0} = \left( \begin{array}{ccccc}
1 & 1 & \frac{1}{2} & \hdots & \frac{1}{(n-1)!}\\
0 & 1 & 1 & \hdots & \frac{1}{(n-2)!}\\
0 & 0 & 1 & \hdots & \frac{1}{(n-3)!}\\
\vdots & \vdots & \vdots & \ddots  & \vdots\\
0 & 0 & 0 & \ldots & 1 \end{array} \right).
\end{align*}
Thus $M_{0}$ has in particular rational entries. Note that we can write $M_{0}=e^{N}$, where $N$ is our notation for the matrix whose non-zero entries are ones on the superdiagonal. In this newly defined basis we have the following theorem.\\ 

\begin{theorem}
\label{InaarF}
The matrix $T$  that transforms functions in the Mellin-Barnes basis $\mathcal{I}$ to the ordered basis $\{f_{n-1}/(2\pi i)^{n-1},f_{n-2}/(2\pi i)^{n-2},\ldots,f_{0}\}$ is given by $T=Q\Phi$. Here $Q$ is the VanderMonde type matrix $Q_{kl} = (k-\frac{n}{2})^{l}/l!$, where $k,l=0,1,\ldots,n-1$, and
\begin{align*}
\Phi =\left( \begin{array}{ccccc}
\phi(0) & \frac{\phi'(0)}{2\pi i} & \frac{\phi''(0)}{2!(2\pi i)^{2}} & \hdots & \frac{\phi^{(n-1)}(0)}{(n-1)!(2\pi i)^{n-1}}\\
0 & \phi(0) & \frac{\phi'(0)}{2\pi i} & \hdots & \frac{\phi^{(n-2)}(0)}{(n-2)!(2\pi i)^{n-2}}\\
0 & 0 & \phi(0) & \hdots & \frac{\phi^{(n-3)}(0)}{(n-3)!(2\pi i)^{n-3}}\\
\vdots & \vdots &  &  & \vdots\\
0 & 0 & 0 & \hdots & \phi(0) \end{array} \right),
\end{align*}
where $\phi$ is the function
\begin{align}
\phi(s) = \frac{\Gamma(\alpha_{1}+s)\cdots\Gamma(\alpha_{n}+s)}{\Gamma(\alpha_{1})\cdots\Gamma(\alpha_{n})}\Gamma(1-s)^{n}.\\
\end{align}
\end{theorem}

\noindent\textbf{Proof. }Let $k\in\{0,1,\ldots,n-1\}$. We see that for $|z|<1$
\begin{align*}
& I_{k}(z) = \frac{(-1)^{n}}{(2\pi i)^{n}}\int_{L}\frac{\Gamma(\alpha_{1}+s)\cdots\Gamma(\alpha_{n}+s)}{\Gamma(\alpha_{1})\cdots\Gamma(\alpha_{n})} \Gamma(-s)^{n} e^{(2k-n)\pi i s} z^{s} ds\\
&= \sum_{m=0}^{\infty} \frac{(-1)^{n}}{(n-1)!} \frac{d^{n-1}}{ds^{n-1}}|_{s=m}\frac{\Gamma(\alpha_{1}+s)\cdots\Gamma(\alpha_{n}+s)}{\Gamma(\alpha_{1})\cdots\Gamma(\alpha_{n})} (s-m)^{n}\Gamma(-s)^{n}\frac{e^{(2k-n)\pi i s}}{(2\pi i)^{n-1}} z^{s}\\
&= \sum_{l=0}^{n-1} \frac{\log^{n-1-l}(z)}{(n-1-l)!}\sum_{m=0}^{\infty} \frac{z^{m}}{l!} \frac{d^{l}}{ds^{l}}|_{s=m} \frac{\Gamma(\alpha_{1}+s)\cdots\Gamma(\alpha_{n}+s)}{\Gamma(\alpha_{1})\cdots\Gamma(\alpha_{n})} (m-s)^{n}\Gamma(-s)^{n}\frac{e^{(2k-n)\pi i s}}{(2\pi i)^{n-1}}\\
&=\sum_{l=0}^{n-1} a_{k,l}(z)\frac{1}{(n-1-l)!}\frac{\log^{n-1-l}(z)}{(2\pi i)^{n-1-l}}
\end{align*}
for suitable analytic functions $a_{k,0},\ldots,a_{k,n-1}$ in a neighborhood of $z=0$ that satisfy in particular
\begin{align*}
a_{k,l}(0) &= \frac{ (2\pi i)^{-l}}{l!} \frac{d^{l}}{ds^{l}}|_{s=0} \frac{\Gamma(\alpha_{1}+s)\cdots\Gamma(\alpha_{n}+s)}{\Gamma(\alpha_{1})\cdots\Gamma(\alpha_{n})}\Gamma(1-s)^{n} e^{(2k-n)\pi i s}\\
&=\sum_{m=0}^{l} \frac{(k-\frac{n}{2})^{m}}{m!}\frac{\phi^{(l-m)}(0)}{(l-m)!(2\pi i)^{l-m}}.
\end{align*}
Here we have used the Leibniz rule. By definition we have $I_{k}(z) = \sum_{l=0}^{n-1} T_{kl} f_{n-1-l}/(2\pi i)^{n-1-l}$ in the Frobenius basis. Since $\log^{k}(z)/k!$ is the only term  in $f_{k}$ which is a power of a logarithm multiplied by a constant term we can apply Proposition \ref{combilogs} to find
\begin{align*}
T_{kl} - a_{k,l}(0) = 0\text{ for }l=0,1,\ldots,n-1.
\end{align*}
\begin{flushright}$\square$\end{flushright}

\begin{proposition}
\label{combilogs}
Let $m\in\mathbb{N}$ and let $a_{0},\ldots,a_{m}$ be analytic functions in a neighborhood of $0$. Suppose that for all $z$ in this neighborhood, with argument in $(0,2\pi)$, we have
\begin{align*}
\sum_{j=0}^{m} a_{j}(z)\log^{j}(z) = 0.
\end{align*}
Then we have $a_{j}(0)=0$ for all $0\leq j\leq m$.
\end{proposition}

\noindent\textbf{Proof. }Suppose the statement of the theorem is untrue. Denote by $0\leq r\leq m$ the largest number such that $a_{r}(0)\neq 0$. We can write
\begin{align*}
a_{r}(z) = -\sum_{j=0}^{r-1} a_{j}(z)\log^{j-r}(z)-\sum_{j=r+1}^{m} a_{j}(z)\log^{j-r}(z).
\end{align*}
Taking the limit $z\to 0$ yields $a_{r}(0)=0$, contradicting our assumption that $r$ was the largest number such that $a_{r}(0)\neq 0$. Here we have used that $\log^{j-r}(z)\to 0$ for $j<r$ and we have used the standard limit $z \log^{j-r}(z)\to 0$ for the terms with $j>r$.  
\begin{flushright}$\square$\end{flushright}

\begin{remark}
\label{combilogremark}
By induction it follows that the analytic functions $a_{j}$ should actually vanish.
\end{remark}

\begin{theorem}
\label{M1stelling}
In the ordered basis $\{f_{n-1}/(2\pi i)^{n-1},f_{n-2}/(2\pi i)^{n-2},\ldots,f_{0}\}$ we have  $M_{1} = \mathbb{I}+u v^T$. Here
\begin{align*}
u =\left( \begin{array}{c} (T^{-1})_{00}\\ (T^{-1})_{10} \\ \vdots \\ (T^{-1})_{(n-1)0} \end{array} \right) \textit{ and }
v = \left( \begin{array}{cccc} \frac{V^{(0)}(0)}{0!}\\ \frac{1}{2\pi i}\frac{V^{(1)}(0)}{1!}\\ \vdots \\ \frac{1}{(2\pi i)^{n-1}}\frac{V^{(n-1)}(0)}{(n-1)!} \end{array} \right)
\end{align*}
and the function $V$ is defined by
\begin{align*}
V(s) = (-1)^{n}\phi(s)e^{-\pi i n s} \prod_{k=1}^{n} (e^{2\pi i s}-e^{-2\pi i \alpha_{k}}).\\
\end{align*}
\end{theorem}

\noindent\textbf{Proof. }From Theorem \ref{MBmonodromy} we obtain in the Mellin-Barnes basis
\begin{align*}
M_{1} &= M_{0}^{-1} M_{\infty}^{-1}\\
&= \left( \begin{array}{ccccc}
(-1)^{n}A_{n} & (-1)^{n}A_{n-1}+\binom{n}{1} & (-1)^{n}A_{n-2}-\binom{n}{2} & \hdots & (-1)^{n}A_{1}\pm\binom{n}{n-1}\\
0 & 1 & 0 & \hdots & 0\\
0 & 0 & 1 & \hdots & 0\\
\vdots & \vdots &  &  & \vdots\\
0 & 0 & 0 & \hdots & 1 \end{array} \right)
\end{align*}
Now we notice that the $(0,l)$th entry of $(M_{1}-\mathbb{I})T$ is
\begin{align*}
((M_{1}-\mathbb{I})T)_{0l} &= \sum_{k=0}^{n-1} (-1)^{n}\left[A_{n-k}-(-1)^{n-k}\binom{n}{k}\right] \frac{(2\pi i)^{n-1-l}}{l!} \phi_{k}^{(l)}(0)
\end{align*}
where
\begin{align*}
\phi_{k}(s) = \frac{\phi(s)}{(2\pi i)^{n-1}}e^{-\pi i n s}e^{2\pi i k s}.
\end{align*}
We see that
\begin{align*}
 \sum_{k=0}^{n-1}\binom{n}{k} (-1)^{n-k}\phi_{k}^{(l)}(0) &= \frac{d^{l}}{ds^{l}}|_{s=0} \frac{\phi(s)}{(2\pi i)^{n-1}}e^{-\pi i n s}\sum_{k=0}^{n-1} \binom{n}{k} (-1)^{n-k}e^{2\pi i k s}\\
&= \frac{d^{l}}{ds^{l}}|_{s=0}\frac{\phi(s)}{(2\pi i)^{n-1}}e^{-\pi i n s}\left((e^{2\pi i s}-1)^{n}-e^{2\pi i n s}\right)\\
&= 0 - \phi_{n}^{(l)}(0) = -A_{0} \phi_{n}^{(l)}(0). 
\end{align*}
Therefore
\begin{align*}
((M_{1}-\mathbb{I})T)_{0l} &=  (2\pi i)^{n-1-l}\frac{(-1)^{n}}{(n-1-l)!} \sum_{k=0}^{n} A_{n-k} \phi_{k}^{(n-1-l)}(0)\\
&= (2\pi i)^{n-1-l}\frac{(-1)^{n}}{l!} \frac{d^{l}}{ds^{l}}|_{s=0}\frac{\phi(s)}{(2\pi i)^{n-1}}e^{-\pi i n s} \prod_{k=1}^{n} (e^{2\pi i s}-e^{-2\pi i \alpha_{k}})\\
&= (2\pi i)^{-l}\frac{V^{(l)}(0)}{l!}.
\end{align*}
Here we used the Leibniz rule. Of course all other entries of $(M_{1}-\mathbb{I})T$ are zero. We conclude that in the ordered basis $\{f_{n-1}/(2\pi i)^{n-1},f_{n-2}/(2\pi i)^{n-2},\ldots,f_{0}\}$ we have
\begin{align*}
M_{1} &= \mathbb{I}+T^{-1}(M_{1}^{\mathcal{I}}-\mathbb{I})T\\
&= \mathbb{I}+\left( \begin{array}{c} (T^{-1})_{00}\\ (T^{-1})_{10} \\ \vdots \\ (T^{-1})_{(n-1)0} \end{array} \right)
 \left( \begin{array}{cccc} \frac{V^{(0)}(0)}{0!}& \frac{1}{2\pi i}\frac{V^{(1)}(0)}{1!}& \hdots & \frac{1}{(2\pi i)^{n-1}}\frac{V^{(n-1)}(0)}{(n-1)!} \end{array} \right).
\end{align*}
Here the superscript $\mathcal{I}$ indicates that the particular matrix is in the Mellin Barnes Basis. 
\begin{flushright}$\square$\end{flushright}

Using Proposition \ref{inverserank1} we get the following corollary. 

\begin{corollary}
\label{M1stelling2}
In the ordered basis $\{f_{n-1}/(2\pi i)^{n-1},f_{n-2}/(2\pi i)^{n-2},\ldots,f_{0}\}$ we have $M_{1}^{\mathcal{F}} = e^{N}+u v^T$. Here
\begin{align*}
u =\left( \begin{array}{c} (T^{-1})_{00}\\ (T^{-1})_{10} \\ \vdots \\ (T^{-1})_{(n-1)0} \end{array} \right) \textit{ and }
v = \left( \begin{array}{cccc} \frac{W^{(0)}(0)}{0!}\\ \frac{1}{2\pi i}\frac{W^{(1)}(0)}{1!}\\ \vdots \\ \frac{1}{(2\pi i)^{n-1}}\frac{W^{(n-1)}(0)}{(n-1)!} \end{array} \right)
\end{align*}
and the function $W$ is defined by $W(s) = (-1)^{n} e^{-2\pi i(\alpha_{1}+\ldots+\alpha_{n})} e^{2\pi i s} V(s)$.
\end{corollary}

\section{The case where $(X-e^{-2\pi i\alpha_{1}})\cdots (X-e^{-2\pi i\alpha_{n}})$ is a product of cyclotomic polynomials}

Theorem \ref{M1stelling} shows us that for large $n$ the expressions for the monodromy matrices seem to become rather cumbersome. Therefore we will, in this chapter, limit our study of the monodromy matrices in the maximally unipotent case to the case where $(X-e^{-2\pi i\alpha_{1}})\cdots (X-e^{-2\pi i\alpha_{n}})$ is a product of cyclotomic polynomials. This is actually not such a big restriction, since it seems to be a case of particular interest (see for example \cite{CYY}). In particular, many Calabi-Yau differential equations are of this form.

\subsection{Polynomials with roots in the cyclotomic field}\index{Cyclotomic field}

\begin{proposition}
\label{polynoomst}
Let $p\in\mathbb{Q}[X]$ be monic and suppose all its roots are roots of unity not equal to $1$. Then there exists a number $r\in\mathbb{N}$ and numbers $a_{1},\ldots,a_{r},b_{1},\ldots,b_{r}\in\mathbb{N}$ such that
\begin{align}
\label{polynoomvorm}
p(X) = \frac{(X^{a_{1}}-1)\cdots (X^{a_{r}}-1)}{(X^{b_{1}}-1)\cdots (X^{b_{r}}-1)}.
\end{align}
\end{proposition}

\noindent\textbf{Proof.} This follows immediately from the fact that the $k^{th}$ cyclotomic polynomial satisfies
\begin{align*}
\phi_{k}(X) = \prod_{d|k} (x^{d}-1)^{\mu(n/d)},
\end{align*}
where $\mu$ denotes the M\"obius function. 
\begin{flushright}$\square$\end{flushright}

\begin{theorem}
\label{gammaproduct}
Let $\alpha_{1},\ldots,\alpha_{n}\in \mathbb{Q}\cap (0,1)$ and suppose that $(X-e^{-2\pi i\alpha_{1}})\cdots (X-e^{-2\pi i\alpha_{n}})$ has integer coefficients. Then there exist a number $r\in\mathbb{N}$ and numbers $a_{1},\ldots,a_{r},b_{1},\ldots,b_{r}\in\mathbb{N}$ such that
\begin{align*}
\prod_{k=1}^{n} \Gamma(\alpha_{k}+s) = C^{-s}\frac{\Gamma(a_{1}s)\cdots \Gamma(a_{r}s)}{\Gamma(b_{1}s)\cdots \Gamma(b_{r}s)} (2\pi)^{\frac{n}{2}}\sqrt{\frac{a_{1}\cdots a_{r}}{b_{1}\cdots b_{r}}}\text{ where }
C = \frac{a_{1}^{a_{1}}\cdots a_{r}^{a_{r}}}{b_{1}^{b_{1}}\cdots b_{r}^{b_{r}}}.
\end{align*}
\end{theorem}

\noindent\textbf{Proof.} By Proposition (\ref{polynoomst}) we find a number $r\in\mathbb{N}$ and numbers $a_{1},\ldots,a_{r},b_{1},\ldots,b_{r}\in\mathbb{N}$ such that
\begin{align*}
\prod_{k=1}^{n} \Gamma(\alpha_{k}+s) = \frac{\left(\prod_{j=0}^{a_{1}-1}\Gamma(\frac{j}{a_{1}}+s)\right)\cdots \left(\prod_{j=0}^{a_{r}-1}\Gamma(\frac{j}{a_{r}}+s)\right)}{\left(\prod_{j=0}^{b_{1}-1}\Gamma(\frac{j}{b_{1}}+s)\right)\cdots \left(\prod_{j=0}^{b_{r}-1}\Gamma(\frac{j}{b_{r}}+s)\right)}.
\end{align*}
This is due to the fact that a bijection can be made between the terms in which the gamma functions are evaluated and the roots of the corresponding polynomials. According to the multiplication theorem for the Gamma function this equals
\begin{align*}
\frac{\left(\Gamma(a_{1}s)(2\pi)^{\frac{a_{1}}{2}} a_{1}^{\frac{1}{2}-a_{1}s}\right)\cdots\left(\Gamma(a_{r}s)(2\pi)^{\frac{a_{r}}{2}} a_{r}^{\frac{1}{2}-a_{r}s}\right)}{\left(\Gamma(b_{1}s)(2\pi)^{\frac{b_{1}}{2}} b_{1}^{\frac{1}{2}-b_{1}s}\right)\cdots \left(\Gamma(b_{r}s)(2\pi)^{\frac{b_{r}}{2}} b_{r}^{\frac{1}{2}-b_{r}s}\right)}\\
=C^{-s}\frac{\Gamma(a_{1}s)\cdots \Gamma(a_{r}s)}{\Gamma(b_{1}s)\cdots \Gamma(b_{r}s)} (2\pi)^{\frac{n}{2}}\sqrt{\frac{a_{1}\cdots a_{r}}{b_{1}\cdots b_{r}}}
\end{align*}
where we have used that $a_{1}+\ldots+a_{r} = n+b_{1}+\ldots+b_{r}$. 
\begin{flushright}$\square$\end{flushright}

\begin{remark}
Notice that we can rewrite this formula as
\begin{align*}
\prod_{k=1}^{n} \Gamma(\alpha_{k}+s) = C^{-s}\frac{\Gamma(a_{1}s+1)\cdots \Gamma(a_{r}s+1)}{\Gamma(b_{1}s+1)\cdots \Gamma(b_{r}s+1)} (2\pi)^{\frac{n}{2}}\sqrt{\frac{b_{1}\cdots b_{r}}{a_{1}\cdots a_{r}}}
\end{align*}
which implies the appealing form
\begin{align}
\label{appeal}
C^{s}\prod_{k=1}^{n} \frac{\Gamma(\alpha_{k}+s)}{\Gamma(\alpha_{k})} = \frac{\Gamma(a_{1}s+1)\cdots \Gamma(a_{r}s+1)}{\Gamma(b_{1}s+1)\cdots \Gamma(b_{r}s+1)}.
\end{align}
\end{remark}

The proof of the following theorem is by Julian Lyczak and Merlijn Staps. 

\begin{proposition}
\label{integerC}
The number $C$ of Theorem \ref{gammaproduct} is an integer. 
\end{proposition}

\noindent\textbf{Proof. }Let $m\in\mathbb{N}$, the number of factors of the product $(X^{b_{1}}-1)\cdots (X^{b_{r}}-1)$ of which $e^{2\pi i/m}$ is a root cannot exceed the number of factors of the product $(X^{a_{1}}-1)\cdots (X^{a_{r}}-1)$ of which $e^{2\pi i/m}$ is a root, otherwise $(X^{a_{1}}-1)\cdots (X^{a_{r}}-1)(X^{b_{1}}-1)^{-1}\cdots (X^{b_{r}}-1)^{-1}$ could not be a polynomial. We conclude that $|\{j:m|a_{j}\}|\geq |\{j:m|b_{j}\}|$ for all $m\in\mathbb{N}$. Now let $p$ be prime and let $k\in\mathbb{N}$. Define $A_{k}=\{a_{j}:p^{k}|a_{j}\}$ and $B_{k}=\{b_{j}:p^{k}|b_{j}\}$ and consider the rational function
\begin{align*}
q(X) = \prod_{a\in A_{k}} (X^{a}-1)/\prod_{b\in B_{k}} (X^{b}-1).
\end{align*}
Suppose $q(X)$ is not a polynomial, then there exists a root of unity $\zeta\neq 1$ such that there are more factors of the form $(X^{b}-1)$ than of the form $(X^{a}-1)$ that have $\zeta$ as a root. This root is of the form $\zeta=e^{2\pi i l/m}$ for some $l,m\in\mathbb{N}$, where $m>1$. In particular, $|\{a\in A_{k}:m|a\}|<|\{b\in B_{k}:m|b\}|$. However, because $|A_{k}|=|\{j:p^{k}|a_{j}\}|\geq |\{j:p^{k}|b_{j}\}|=|B_{k}|$ we must have $gcd(m,p)=1$, and this would imply $|\{j:p^{k}m|a_{j}\}|<|\{j:p^{k}m|b_{j}\}|$, which is a contradiction. We must conclude that $q(X)$ is a polynomial, thus by comparing degrees we have
\begin{align*}
\sum_{a\in A_{k}} a \geq \sum_{b\in B_{k}} b.
\end{align*}
Denote by $\mathcal{A}_{j}$ the largest integer such that $p^{\mathcal{A}_{j}}| a_{j}$ and by $\mathcal{B}_{j}$ the largest integer such that $p^{\mathcal{B}_{j}} | b_{j}$. The theorem is now proved by the observation that
\begin{align*}
\sum_{j=1}^{r} \mathcal{A}_{j} a_{j} = \sum_{k=1}^{\infty} \sum_{a\in A_{k}} a\geq \sum_{k=1}^{\infty} \sum_{b\in B_{k}} b = \sum_{j=1}^{r} \mathcal{B}_{j} b_{j}.
\end{align*}
\begin{flushright}$\square$\end{flushright}

\begin{corollary}
\label{faculteit}
Let $r\in\mathbb{N}$ and let $a_{1},\ldots,a_{n},b_{1},\ldots,b_{n}\in \mathbb{N}$. Suppose that
\begin{align}
\label{poly1}
\frac{(X^{a_{1}}-1)\cdots (X^{a_{r}}-1)}{(X^{b_{1}}-1)\cdots (X^{b_{r}}-1)}
\end{align}
is a polynomial. Then $\frac{a_{1}!\cdots a_{r}!}{b_{1}!\cdots b_{r}!}$ is an integer. 
\end{corollary}

\noindent\textbf{Proof. }Notice that by multiplying with $(X-1)$ we may assume (\ref{poly1}) to be non-constant. Without loss of generality (\ref{poly1}) is irreducible (this follows from Proposition \ref{polynoomst}). Thus there exists a $\mathcal{N}\in\mathbb{N}$ such that $\{\alpha_{1},\ldots,\alpha_{n}\}=\{m/\mathcal{N}:0<m<\mathcal{N},gcd(m,\mathcal{N})=1\}$. It follows from (\ref{appeal}) that
\begin{align*}
\frac{a_{1}!\cdots a_{r}!}{b_{1}!\cdots b_{r}!} = \alpha_{1}\cdots \alpha_{n} \frac{a_{1}^{a_{1}}\cdots a_{r}^{a_{r}}}{b_{1}^{b_{1}}\cdots b_{r}^{b_{r}}}.
\end{align*}
Let $p$ be a prime divisor of $\mathcal{N}$ and denote by $m$ its multiplicity. We follow the proof of Proposition \ref{integerC} untill we define the polynomial
\begin{align*}
q(X) = \prod_{a\in A_{k}} (X^{a}-1)/\prod_{b\in B_{k}} (X^{b}-1).
\end{align*}
for $k\leq m$ (with same notation). Notice that indeed there must exist an $a_{j}$ such that $p^{m}|a_{j}$ because $e^{2\pi i/\mathcal{N}}$ must be a root of our original polynomial. In this case, we can reason that $e^{-2\pi i\alpha_{j}}$ must be a root of $q(X)$, this is because it is a root of our original polynomial and cannot be a root of any factor not corresponding to $A_{k}$. By comparing degrees we conclude that
\begin{align*}
\sum_{a\in A_{k}} a \geq n+\sum_{b\in B_{k}} b.
\end{align*}
We obtain
\begin{align*}
-m n+\sum_{j=1}^{r} \alpha_{j} a_{j} = -m n+\sum_{k=1}^{\infty} \sum_{a\in A_{k}} a\geq \sum_{k=1}^{\infty} \sum_{b\in B_{k}} b = \sum_{j=1}^{r} \beta_{j} b_{j}
\end{align*}
which proves our corollary. 
\begin{flushright}$\square$\end{flushright}

\subsection{A general expression for the monodromy matrices of the maximally unipotent case}

\noindent 
If we would instead of the generalized hypergeometric equation have considered the equation
\begin{align}
\label{hypergeoC}
\theta^{n}f=Cz(\theta+\alpha_{1})\cdots (\theta+\alpha_{n})f
\end{align}
then a solution $f$ to this equation for $C=1$, i.e. of the hypergeometric case, induces the solution $f(C z)$ for general $C\in\mathbb{C}\setminus\{0\}$. In other words, normalization of $z$ provides us with solutions to a related differential equation. Let us use our knowledge of the hypergeometric equation to find `a Frobenius basis' for (\ref{hypergeoC}). Denote this Frobenius basis by $f_{0}^{C},\ldots,f_{n-1}^{C}$. We know that a basis of solutions is given by $f_{0}(Cz),\ldots,f_{n-1}(Cz)$. Notice that
\begin{align*}
f_{j}(Cz) &= \frac{\log^{j}(Cz)}{j!}+\sum_{m=0}^{j} \frac{\log^{m}(Cz)}{m!}h_{m}(Cz)\\
&= \sum_{m=0}^{j} \frac{\log^{m}(z)}{m!} \frac{\log^{j-m}(C)}{(j-m)!}+\sum_{m=0}^{j} \frac{(\log(z)+\log(C))^{m}}{m!}h_{m}(Cz)\\
&= \sum_{m=0}^{j} \frac{\log^{j-m}(C)}{(j-m)!} f_{j}^{C}(z).
\end{align*}
We conclude that
\begin{align*}
\left(\begin{array}{c} f_{n-1}^{C}(z)/(2\pi i)^{n-1} \\ \vdots \\ f_{0}^{C}(z) \end{array}\right) = C^{-\frac{N}{2\pi i}} \left(\begin{array}{c} f_{n-1}(Cz)/(2\pi i)^{n-1} \\ \vdots \\ f_{0}(Cz) \end{array}\right).
\end{align*}
Again $N$ is the matrix who's only nonzero components are ones on the superdiagonal. Notice that in this case our monodromy group is generated by $M_{0},M_{1/C}$ and $M_{\infty}$. \\

From now on we choose $C$ to be the constant from the previous paragraph, that is
\begin{align*}
C = \frac{a_{1}^{a_{1}}\cdots a_{n}^{a_{n}}}{b_{1}^{b_{1}}\cdots b_{n}^{b_{n}}}.
\end{align*}

\begin{theorem}
Let $\alpha_{1},\ldots,\alpha_{n}\in \mathbb{Q}\cap (0,1)$ and suppose that $(X-e^{-2\pi i\alpha_{1}})\cdots (X-e^{-2\pi i\alpha_{n}})$ has integer coefficients. Then the solution $f_{0}^{C}$ of (\ref{hypergeoC}) has integer coefficients in its powerseries expansion. 
\end{theorem}

\noindent\textbf{Proof. }From the above discussion we infer that
\begin{align*}
f_{0}^{C}(z) = {_{n}F_{n-1}}(\alpha_{1},\ldots,\alpha_{n};1,\ldots,1|Cz) = \sum_{m=0}^{\infty} \frac{(a_{1}m)!\cdots (a_{r}m)!}{(b_{1}m)!\cdots (b_{r}m)!}\frac{z^{m}}{m!^{n}},
\end{align*}
where we have used (\ref{appeal}). Without loss of generality $(X-e^{-2\pi i\alpha_{1}})\cdots (X-e^{-2\pi i\alpha_{n}})$ is irreducible. Let $p\leq m$ be prime. Let $\mathcal{N}$ be as in corollary \ref{faculteit}. Suppose $p\not|\mathcal{N}$. We have
\begin{align*}
\frac{(a_{1}m)!\cdots (a_{r}m)!}{(b_{1}m)!\cdots (b_{r}m)!} &= \left(\prod_{k=1}^{n} \prod_{l=0}^{m-1} \frac{\alpha_{k}+l}{m}\right) \frac{(a_{1}m)^{a_{1}m}\cdots (a_{r}m)^{a_{r}m}}{(b_{1}m)^{b_{1}m}\cdots (b_{r}m)^{b_{r}m}}\\
&= \left(\prod_{k=1}^{n} \prod_{l=0}^{m-1}  \frac{\mathcal{N}\alpha_{k}+\mathcal{N} l}{\mathcal{N}}\right) \left(\frac{a_{1}^{a_{1}}\cdots a_{r}^{a_{r}}}{b_{1}^{b_{1}}\cdots b_{r}^{b_{r}}}\right)^{m}.
\end{align*}
Because $\text{gcd}(p,\mathcal{N})=1$ we have $\{0,\mathcal{N},2\mathcal{N},\ldots,(p^{l}-1)\mathcal{N}\}\equiv\{0,1,\ldots,p^{l}-1\} \mod{p^{l}}$. Thus at least $[m/p^{l}]$ of $\mathcal{N} \alpha_{k}, \mathcal{N} \alpha_{k}+\mathcal{N},\ldots,\mathcal{N} \alpha_{k}+(m-1)\mathcal{N}$ must be divisible by $p^{l}$. We conclude that
\begin{align*}
p^{n([m/p]+[m/p^{2}]+\ldots)}|\prod_{k=1}^{n} \prod_{l=0}^{m-1}  (\mathcal{N}\alpha_{k}+\mathcal{N} l),
\end{align*}
and this is enough. Now suppose $p|\mathcal{N}$ with multiplicity $e$. We notice that
\begin{align*}
\frac{(a_{1}m)!\cdots (a_{r}m)!}{(b_{1}m)!\cdots (b_{r}m)!} &=  \left(\prod_{k=1}^{n} \prod_{l=0}^{m-1}  \frac{\mathcal{N}\alpha_{k}+\mathcal{N} l}{\mathcal{N}\alpha_{k}}\right) \left(\frac{a_{1}!\cdots a_{r}!}{b_{1}!\cdots b_{r}!}\right)^{m}
\end{align*}
We should prove that
\begin{align*}
p|\frac{a_{1}!\cdots a_{r}!}{b_{1}!\cdots b_{r}!}.
\end{align*}
If this is not the case then we deduce from the proof of corollary \ref{faculteit} that
\begin{align*}
(X-e^{-2\pi i\alpha_{1}})\cdots (X-e^{-2\pi i\alpha_{n}}) = \prod_{a\in A_{e}} (X^{a}-1)/\prod_{b\in B_{e}} (X^{b}-1).
\end{align*}
Thus $(X-e^{-2\pi i\alpha_{1}})\cdots (X-e^{-2\pi i\alpha_{n}})=q(X^{p^{e}})$ for some polynomial $q$ that must necessarily be cyclotomic and irreducible. We conclude that there must exist an $\mathcal{M}\in\mathbb{N}$ such that $\varphi(\mathcal{N})=n\varphi(\mathcal{M})$, where $\varphi$ is the Euler totient function\index{Euler's totient function}. Also we deduce that $p^{e}|n$. Since $e^{2\pi i p^{e}/\mathcal{N}}$ is a root of $q$ we must have $\mathcal{N}/p^{e}|\mathcal{M}$. Hence
\begin{align*}
\varphi(\mathcal{N}) = n\varphi(\mathcal{M}) \geq n\varphi(\mathcal{N}/p^{e}) = \varphi(\mathcal{N}) \frac{n}{p^{e}} \frac{p}{p-1} > \varphi(\mathcal{N}),
\end{align*}
a contradiction. 
\begin{flushright}$\square$\end{flushright}

The authors of \cite{CYY} point out that this result holds for all Picard-Fuchs equations (i.e. the $n=4$ case), it is actually used as part of the definition of a Calabi-Yau type differential equation by the authors of \cite{AESZ}. A folklore conjecture that goes back to Bombieri and Dwork states that all power series $y_{0}(z)\in\mathbb{Z}[[z]]$ that satisfy a homogeneous linear differential equation have a geometrical origin.\\

Matrices that have the form of $\Phi$ from Theorem \ref{InaarF} have a certain homomorphism\index{Homomorphism} property. Explicitly, for a function $C(s)$ we have
\begin{align*}
& \left( \begin{array}{ccccc}
\phi(0) & \frac{\phi'(0)}{2\pi i} & \frac{\phi''(0)}{2!(2\pi i)^{2}} & \hdots & \frac{\phi^{(n-1)}(0)}{(n-1)!(2\pi i)^{n-1}}\\
0 & \phi(0) & \frac{\phi'(0)}{2\pi i} & \hdots & \frac{\phi^{(n-2)}(0)}{(n-2)!(2\pi i)^{n-2}}\\
0 & 0 & \phi(0) & \hdots & \frac{\phi^{(n-3)}(0)}{(n-3)!(2\pi i)^{n-3}}\\
\vdots & \vdots &  &  & \vdots\\
0 & 0 & 0 & \hdots & \phi(0) \end{array} \right)
\left( \begin{array}{ccccc}
C(0) & \frac{C'(0)}{2\pi i} & \frac{C''(0)}{2!(2\pi i)^{2}} & \hdots & \frac{C^{(n-1)}(0)}{(n-1)!(2\pi i)^{n-1}}\\
0 & C(0) & \frac{C'(0)}{2\pi i} & \hdots & \frac{C^{(n-2)}(0)}{(n-2)!(2\pi i)^{n-2}}\\
0 & 0 & C(0) & \hdots & \frac{C^{(n-3)}(0)}{(n-3)!(2\pi i)^{n-3}}\\
\vdots & \vdots &  &  & \vdots\\
0 & 0 & 0 & \hdots & C(0) \end{array} \right)\\
&= \left( \begin{array}{ccccc}
\phi_{C}(0) & \frac{\phi_{C}'(0)}{2\pi i} & \frac{\phi_{C}''(0)}{2!(2\pi i)^{2}} & \hdots & \frac{\phi_{C}^{(n-1)}(0)}{(n-1)!(2\pi i)^{n-1}}\\
0 & \phi_{C}(0) & \frac{\phi_{C}'(0)}{2\pi i} & \hdots & \frac{\phi_{C}^{(n-2)}(0)}{(n-2)!(2\pi i)^{n-2}}\\
0 & 0 & \phi_{C}(0) & \hdots & \frac{\phi_{C}^{(n-3)}(0)}{(n-3)!(2\pi i)^{n-3}}\\
\vdots & \vdots &  &  & \vdots\\
0 & 0 & 0 & \hdots & \phi_{C}(0) \end{array} \right)
\end{align*}
where $\phi_{C}(s) = \phi(s) C(s)$. Notice that the second matrix in the product is simply $C^{\frac{N}{2\pi i}}$ when $C(s)$ is defined to be $C^{s}$. The results we have found so far adapt naturally to the new basis (where $z$ is normalized with $C$), we simply substitute $\phi$ by $\phi_{C}$ (compare this with Theorem \ref{InaarF}). It should be clear why this basis is interesting, with our particular choice of $C$ we have the appealing form
\begin{align*}
\phi_{C}(s) =  \frac{\Gamma(a_{1}s+1)\cdots \Gamma(a_{r}s+1)}{\Gamma(b_{1}s+1)\cdots \Gamma(b_{r}s+1)}\Gamma(1-s)^{n}.
\end{align*}

\begin{definition}
Let $j\in\mathbb{N}$. By $\pi_{j}$ we denote the set of integer partitions\index{Integer partition} of $j$, i.e. the set of finite (not necessarily strictly) decreasing sequences of natural numbers $p_{1},p_{2},\ldots$ such that $p_{1}+p_{2}+\ldots=j$. Any function $g$ whose domain contains $\mathbb{N}$ can be extended to partitions by multiplication, i.e. $g(p) = g(p_{1})g(p_{2})\cdots$. Additionally, we define $\pi_{0}=\{0\}$ and $g(0)=1$. 
\end{definition}

The following theorem will provide us with a practical method to obtain the monodromy matrices in the ordered basis $f_{n-1}^{C}/(2\pi i)^{n-1},\ldots, f_{1}^{C}/(2\pi i), f_{0}^{C}$.\\

\begin{theorem}
\label{main?}
\textbf{(Main Theorem)}\\
Let $\alpha_{1},\ldots,\alpha_{n}\in \mathbb{Q}\cap (0,1)$ and suppose that $(X-e^{-2\pi i\alpha_{1}})\cdots (X-e^{-2\pi i\alpha_{n}})$ is a product of cyclotomic polynomials. Let $r\in\mathbb{N}$ and $a_{1},\ldots,a_{r},b_{1},\ldots,b_{r}\in\mathbb{N}$ be as in Theorem \ref{gammaproduct} and define $\zeta(1)=0$ for convenience. In the ordered basis $f_{n-1}^{C}/(2\pi i)^{n-1},\ldots,f_{1}^{C}/(2\pi i),f_{0}^{C}$ of (\ref{hypergeoC}) we have $M_{1/C} = \mathbb{I} - v_{-} v_{+}^{T}$, where
\begin{align*}
\label{entries}
v_{-,j} &= \sum_{l=0}^{n-1-j} c_{l+j} \sum_{p\in\pi_{l}} \frac{1}{M(p)} c_{p}^{-}\frac{\zeta(p)}{(2\pi i)^{p}}\text{ and }v_{+,j} = \sum_{p\in \pi_{j}} \frac{1}{M(p)} c_{p}^{+}\frac{\zeta(p)}{(2\pi i)^{p}}
\end{align*}
for $j=0,1,\ldots,n-1$. Here the coefficients $c_{j},c_{j}^{\pm}\in\mathbb{Q}$ are given by $c_{0}^{\pm}=1$ and
\begin{align*}
c_{j}^{\pm} &=  \frac{1}{j}\left(\pm n- (\pm 1)^{j}\sum_{m=1}^{r} (a_{m}^{j}-b_{m}^{j})\right)\text{ and }c_{j} = \frac{1}{(n-1)!}\frac{a_{1}\cdots a_{r}}{b_{1}\cdots b_{r}}\frac{d^{j}}{dz^{j}}\left. \prod_{m=1}^{n-1} \left(z-m+\frac{n}{2}\right)\right|_{z=0} 
\end{align*}
(the definition for $c_{j}$ also being valid for $j=0$) and the function $M:\pi_{0}\cup \pi_{1} \cup\cdots\to\mathbb{N}$ by $M(p_{1},p_{2},\cdots) = |\{k:p_{k}=1\}|! |\{k:p_{k}=2\}|!\cdots$\\

\noindent In particular, all matrices in the corresponding monodromy group have their entries in $\mathbb{Q}(\zeta(3) (2\pi i)^{-3},\zeta(5) (2\pi i)^{-5},\ldots,\zeta(m)(2\pi i)^{-m})$, with $m$ the largest odd number below $n$.\\
\end{theorem}

\noindent\textbf{Proof. }We use the function $V$ from theorem (\ref{M1stelling}). After conjugation with the matrix $C^{\frac{N}{2\pi i}}$ we have the same theorem but with function $\phi_{C}(s)=C^{s}\phi(s)$ instead. Notice that
\begin{align*}
(-1)^{n}e^{\pi i(\alpha_{1}+\ldots+\alpha_{n})} V_{C}(s) &:= \phi(s) C^{s}\prod_{k=1}^{n} (e^{\pi i(\alpha_{k}+s)}-e^{-\pi i(\alpha_{k}+s)})\\
&= (2\pi i)^{n}\Gamma(1-s)^{n}C^{s}\prod_{k=1}^{n} \frac{1}{\Gamma(\alpha_{k})\Gamma(1-\alpha_{k}-s)}\\
&= (2\pi i)^{n} \frac{\Gamma(1-s)^{n}}{\Gamma(\alpha_{1})^{2}\cdots \Gamma(\alpha_{n})^{2}} \frac{\Gamma(1-b_{1}s)\cdots \Gamma(1-b_{r}s)}{\Gamma(1-a_{1}s)\cdots \Gamma(1-a_{r}s)}\\
&= i^{n} \Gamma(1-s)^{n} \frac{a_{1}\cdots a_{r}}{b_{1}\cdots b_{r}}\frac{\Gamma(1-b_{1}s)\cdots \Gamma(1-b_{r}s)}{\Gamma(1-a_{1}s)\cdots \Gamma(1-a_{r}s)}
\end{align*}
We remark that one must have $\alpha_{1}+\ldots+\alpha_{n}=\frac{n}{2}$. Using the formula
\begin{align*}
\log \Gamma(1+s) &= -\gamma s+\sum_{p=2}^{\infty} \frac{(-1)^{p}}{p} \zeta(p) s^{p}
\end{align*}
yields
\begin{align*}
(-1)^{n}\frac{b_{1}\cdots b_{r}}{a_{1}\cdots a_{r}} V_{C}(s) &= \exp\left(\sum_{p=2}^{\infty} c_{p}^{+} \zeta(p) s^{p}\right)\\
&= 1+\sum_{r=1}^{\infty} \frac{1}{r!} \left(\sum_{p_{1}=1}^{\infty} c_{p_{1}}^{+} \zeta(p_{1}) s^{p_{1}}\right) \left(\sum_{p_{2}=1}^{\infty} c_{p_{2}}^{+} \zeta(p_{2}) s^{p_{2}}\right) \cdots \left(\sum_{p_{r}=1}^{\infty} c_{p_{r}}^{+} \zeta(p_{r}) s^{p_{r}}\right)\\
&= 1+\sum_{j=1}^{\infty} s^{j} \sum_{r=1}^{j} \frac{1}{r!} \sum_{p_{1}+\cdots+p_{r}=j}  c_{p_{1}}^{+} \zeta(p_{1})\cdots c_{p_{r}}^{+} \zeta(p_{r})\\
&= \sum_{j=0}^{\infty} \left(\sum_{p\in\pi_{j}} \frac{1}{M(p)} c_{p}^{+}\frac{\zeta(p)}{(2\pi i)^{p}}\right) (2\pi i s)^{j},
\end{align*}
where $p c_{p}^{+} = n-(a_{1}^{p}+\ldots+a_{r}^{p}-b_{1}^{p}-\ldots-b_{r}^{p})$. To complete the proof we will have to know the inverse of $Q\Phi C^{\frac{N}{2\pi i}}$. The inverse of $\Phi C^{\frac{N}{2\pi i}}$ is obvious from the homomorphism property of this type of matrix. We remark that the inverse of $Q$ is determined by
\begin{align*}
\prod_{m=0,m\neq k}^{n-1}\frac{(z-m+\frac{n}{2})}{k-j} = \frac{(Q^{-1})_{0,k}}{0!}+\frac{(Q^{-1})_{1,k}}{1!} z+\ldots+\frac{(Q^{-1})_{n-1,k}}{(n-1)!} z^{n-1}.
\end{align*}
Fortunately we will only need the first column. We find
\begin{align*}
(Q^{-1})_{l,0} = \frac{(-1)^{n-1}}{(n-1)!}\frac{d^{l}}{dz^{l}}|_{z=0} \prod_{m=1}^{n-1} (z-m+\frac{n}{2}).
\end{align*}
We notice that
\begin{align*}
\frac{1}{\phi_{C}(s)} &= \Gamma(1-s)^{-n} \frac{\Gamma(b_{1}s)\cdots\Gamma(b_{r}s)}{\Gamma(a_{1}s)\cdots\Gamma(a_{r}s)}\\
&= \exp\left(\sum_{p=2} c_{p}^{-}\zeta(p) s^{p}\right)\\
&= \sum_{j=0}^{\infty} \left(\sum_{p\in\pi_{j}} \frac{1}{M(p)} c_{p}^{-}\frac{\zeta(p)}{(2\pi i)^{p}}\right) (2\pi i s)^{j},
\end{align*}
where $p c_{p}^{-} = -n-(-1)^{p}(a_{1}^{p}+\ldots+a_{r}^{p}-b_{1}^{p}-\ldots-b_{r}^{p})$. It follows that
\begin{align*}
(n-1)!(-1)^{n-1}(C^{-\frac{N}{2\pi i}}\Phi^{-1}Q^{-1})_{j,0} &=\sum_{l=0}^{n-1-j} \frac{d^{l+j}}{dz^{l+j}}|_{z=0} \prod_{m=1}^{n-1} \left(z-m+\frac{n}{2}\right) \sum_{p\in\pi_{j}} \frac{1}{M(p)} c_{p}^{-}\frac{\zeta(p)}{(2\pi i)^{p}}
\end{align*}
The last part of the theorem follows from the fact that $M_{0}$ has integer coefficients and 
\begin{align*}
\frac{\zeta(2p)}{(2\pi i)^{2p}} = -\frac{B_{2p}}{2(2p)!}
\end{align*}
where $B_{2p}$ is the $2p$-th Bernoulli number\index{Bernoulli number}.
\begin{flushright}$\square$\end{flushright}

\begin{remark}
Notice that the above theorem produces a practical method to determine monodromy matrices. Given $\alpha_{1},\ldots,\alpha_{n}\in\mathbb{Q}\cap (0,1)$ one has to write the corresponding polynomial in the form (\ref{polynoomvorm}) and then simply calculate the coefficients $c_{j}^{\pm},c_{j}$. 
\end{remark}

\begin{remark}
For the last part of the main theorem the $\alpha_{k}$ need not actually lie in $(0,1)$ as can be seen from the multiplicative property of the gamma function and the homomorphism property of the $\Phi$ matrix. 
\end{remark}
We point out that  in the Frobenius basis $f_{0},f_{1},\ldots,f_{n-1}$ the monodromy matrices can be obtained by a trivial transformation, namely inverting the conjugation by $C^{\frac{N}{2\pi i}}$. Hence the entries are in $\mathbb{Q}(\log(C)(2 \pi i)^{-1},\zeta(3) (2\pi i)^{-3},\zeta(5) (2\pi i)^{-5}\ldots,\zeta(m)(2\pi i)^{-m})$, with $m$ the largest odd number below $n$.

\subsection{Applications of the main theorem}

As one can check the case $n=2$ yields
\begin{align*}
M_{1} = \left(\begin{array}{cc} 1 & 0\\ -\frac{a_{1}\cdots a_{r}}{b_{1}\cdots b_{r}} & 1\end{array}\right) . 
\end{align*}
The results are summarized in the following table. 

\begin{align*}
\begin{array}{ l || c }
  \text{Case} & \frac{a_{1}\cdots a_{r}}{b_{1}\cdots b_{r}} \\
\hline
  (z+1)^{2} = \frac{(z^{2}-1)^{2}}{(z-1)^{2}} & 4 \\
 (z^{2}+z+1) = \frac{z^{3}-1}{z-1} & 3 \\
 (z^{2}+1) = \frac{z^{4}-1}{z^{2}-1} & 2 \\
 (z^{2}-z+1) = \frac{(z^{6}-1)(z-1)}{(z^{3}-1)(z^{2}-1)} & 1 
\end{array}
\end{align*}

\noindent Let us look at the case $n=3$. Using the identity $c_{2}^{-}+3=c_{2}^{+}$ we obtain the matrix
\begin{align*}
M_{1} = \left(\begin{array}{ccc} 1+b d & 0 & -b^{2} d\\ 0 & 1 & 0\\ -d & 0 &1+ b d\end{array}\right)
\end{align*}
where 
\begin{align*}
b = \frac{c_{2}^{+}}{24}\text{ and }d=\frac{a_{1}\cdots a_{r}}{b_{1}\cdots b_{r}}. 
\end{align*}
 All the corresponding cases are worked out in the following table.\\

\begin{align*}
\begin{array}{ l || c | c | c }
  \text{Case} & C & 24 b & d/2 \\
\hline
  (z+1)^{3} = \frac{(z^{2}-1)^{3}}{(z-1)^{3}} & 64 & -3 & 4\\
 (z^{2}+z+1)(z+1) = \frac{(x^{2}-1)(x^{3}-1)}{(x-1)^{2}} & 108 & -4 & 3\\
 (z^{2}+1)(z+1) = \frac{x^{4}-1}{x-1} & 256 & -6 & 2\\
 (z^{2}-z+1)(z+1) = \frac{x^{6}-1}{x^{3}-1} & 1728 & -12 & 1
\end{array}
\end{align*}
From this table we read off that $bd=-1$ in all cases and we deduce the even nicer form
\begin{align*}
M_{1} = \left(\begin{array}{ccc} 0 & 0 & -1/d \\ 0 & 1 & 0\\ -d & 0 & 0\end{array}\right).
\end{align*}

\noindent Let us apply the theorem to the case $n=4$. This case corresponds to the Picard-Fuchs equation, given by
\begin{align}
[\theta^{4}-C z(\theta-A)(\theta+A-1)(\theta-B)(\theta+B-1)]f = 0.
\end{align}
These differential equations arise from Calabi-Yau threefolds\index{Calabi-Yau threefold} (see \cite{CYY}). Let us apply the main theorem, using that $c_{2}^{-}+4 = c_{2}^{+}$ and $c_{3}^{-} = -c_{3}^{+}$ we can write $M_{1}$ as
\begin{align*}
\left(\begin{array}{cccc} 1+a & 0 & ab/d & a^{2}/d\\
-b & 1 & -b^{2}/d & -ab/d\\
0 & 0 & 1 & 0\\
-d & 0 & -b & 1-a \end{array}\right)
\end{align*}
when we identify
\begin{align*}
d &= \frac{a_{1}\cdots a_{r}}{b_{1}\cdots b_{r}}, a = d c_{3}^{+}\frac{\zeta(3)}{(2\pi i)^{3}}\text{ and }b = -\frac{d c_{2}^{+}}{24}.
\end{align*}

The authors of \cite{CYY} point out that the entries of $M_{1/C}$ contain geometric invariants belonging to the corresponding Calabi-You threefolds. The 14 corresponding cases are worked out in the following table. 

\begin{align*}
\begin{array}{ l | r || c | c | c | c }
  \text{Case} & \text{Polynomial} & C & d & 24 b & (2\pi i)^{3} a/\zeta(3) \\
\hline
(1/5,2/5,3/5,4/5) & \frac{X^{5}-1}{X-1} & 3025 & 5 & 50 & -200\\
(1/10,3/10,7/10,9/10) & \frac{(X-1)(X^{10}-1)}{(X^{2}-1)(X^{5}-1)} & 800 000 & 1 & 34 & -288\\
(1/2,1/2,1/2,1/2) & \frac{(X^{2}-1)^{4}}{(X-1)^{4}} & 256 & 16 & 64 & -128\\
(1/3,1/3,2/3,2/3) & \frac{(X^{3}-1)^{2}}{(X-1)^{2}} & 729 & 9 & 54 & -144\\
(1/3,1/2,1/2,2/3) & \frac{(X^{2}-1)^{2}(X^{3}-1)}{(X-1)^{3}} & 432 & 12 & 60 & -144\\
(1/4,1/2,1/2,3/4) & \frac{(X^{2}-1)(X^{4}-1)}{(X-1)^{2}} & 1024 & 8 & 56 & -176\\
(1/8,3/8,5/8,7/8) & \frac{X^{8}-1}{X^{4}-1} & 65 536 & 2 & 44 & -296\\
(1/6,1/3,2/3,5/6) & \frac{X^{6}-1}{X^{2}-1} & 11664 & 3 & 42 & -204\\
(1/12,5/12,7/12,11/12) & \frac{(X^{2}-1)(X^{12}-1)}{(X^{4}-1)(X^{6}-1)} & 2 985 984 & 1 & 46 & -484\\
(1/4,1/4,3/4,3/4) & \frac{(X^{4}-1)^{2}}{(X^{2}-1)^{2}} & 496 & 4 & 40 & -144\\
(1/4,1/3,2/3,3/4) & \frac{(X^{3}-1)(X^{4}-1)}{(X-1)(X^{2}-1)} & 1728 & 6 & 48 & -156\\
(1/6,1/4,3/4,5/6) & \frac{(X-1)(X^{4}-1)(X^{6}-1)}{(X^{2}-1)^{2} (X^{3}-1)} & 27 648 & 2 & 32 & -156\\
(1/6,1/6,5/6,5/6) & \frac{(X-1)^{2} (X^{6}-1)^{2}}{(X^{2}-1)^{2} (X^{3}-1)^{2}} & 186 624 & 1 & 22 & -120\\
(1/6,1/2,1/2,5/6) & \frac{(X^{2}-1)(X^{6}-1)}{(X-1)(X^{3}-1)} & 6912 & 4 & 52 & -256
\end{array}
\end{align*}
This is in agreement with the results of \cite{CYY}. 

\pagebreak

\phantomsection
\cleardoublepage

\end{document}